\documentclass{IEE-CTA}

\usepackage{epstopdf}
\usepackage{amsmath}
\usepackage{amssymb}
\usepackage{mdwmath}
\usepackage{mathrsfs}
\usepackage{subfig}
\usepackage[Gray,squaren,thinqspace,thinspace]{SIunits}
\usepackage{pgfplots}
\usepackage{tikz}
\usepackage{multirow}

\newtheorem{theorem}{Theorem}{}
{}
{}

\newcommand{\probspac}{\mathscr{P}}
\newcommand{\parspac}{\mathcal{P}}
\newcommand{\oparspac}{\overline{\parspac}_\infty}
\newcommand{\oparspactwo}{\overline{\parspac}_2}

\newcommand{\MVP}{\vec{A}}
\newcommand{\pp}{\textbf{P}}
\newcommand{\Npar}{N_\mathrm{par}}
\newcommand{\pone}{P_1}
\newcommand{\ptwo}{P_2}
\newcommand{\pthree}{P_3}
\newcommand{\ppexpv}{\overline{\pp}}
\newcommand{\pexpv}{\overline{P}}
\newcommand{\phat}{\hat{\pp}}
\newcommand{\ppup}{P^{\mathrm{u}}}
\newcommand{\pplow}{P^{\mathrm{l}}}
\newcommand{\EE}{\mathbb{E}}
\newcommand{\MM}{\mathbb{M}}
\newcommand{\std}{\mathrm{std}}
\newcommand{\rDelta}{\Delta}
\newcommand{\ddelta}{\boldsymbol{\Delta}}
\newcommand{\dlow}{\rDelta^{\mathrm{l}}}
\newcommand{\dup}{\rDelta^{\mathrm{u}}}
\newcommand{\ddlow}{\ddelta^{\mathrm{l}}}
\newcommand{\ddup}{\ddelta^{\mathrm{u}}}
\newcommand{\dbound}{\rDelta^{\mathrm{b}}}

\begin{document}

\title{Robust Optimization Formulations for the Design of an Electric Machine}

\author{\au{Z. Bontinck$^{1,2,*}$}
\au{O. Lass$^3$}
\au{S. Sch\"ops$^{1,2}$}
\au{H. De Gersem$^{2}$}
\au{S. Ulbrich$^{3}$}
\au{O. Rain$^{4}$}
}

\address{\add{1}{Graduate School of Computational Engineering, Technische Universit\"at Darmstadt, Dolivostra\ss e 15, 64293 Darmstadt, Germany}
\add{2}{Institut f\"ur Theorie Elektromagnetischer Felder, Technische Universit\"at Darmstadt, Schlo\ss gartenstra\ss e 8, 64289 Darmstadt, Germany}
\add{3}{Chair of Nonlinear Optimization, Department of Mathematics, Technische Universit\"at Darmstadt, Dolivostra\ss e 15, 64293 Darmstadt, Germany}
\add{4}{Robert Bosch GmbH, 70049 Stuttgart, Germany}
\email{bontinck@gsc.tu-darmstadt.de}}

\begin{abstract}
In this paper two formulations for the robust optimization of the size of the permanent magnet in a synchronous machine are discussed. The optimization is constrained by a partial differential equation to describe the electromagnetic behavior of the machine. The need for a robust optimization procedure originates from the fact that optimization parameters have deviations. The first approach, i.e., \textcolor{red}{worst-case} optimization, makes use of local sensitivities. The second approach takes into account expectation values and standard deviations. The latter are associated with global sensitivities.  The geometry parametrization is elegantly handled thanks to the introduction of an affine decomposition. Since the stochastic quantities are determined by tools from uncertainty quantification (UQ) and thus require a lot of finite element evaluations, model order reduction is used in order to increase the efficiency of the procedure. It is shown that both approaches are equivalent if a linearization is carried out. \textcolor{black}{This finding is supported by the application on an electric machine. The optimization algorithms used are sequential quadratic programming, particle swarm optimization and genetic algorithm}. While both formulations reduce the size of the magnets, the UQ based optimization approach is less pessimistic with respect to deviations and yields smaller magnets.
\end{abstract}

\maketitle

\section{Introduction}
Electric machines are often subjected to optimization in order to improve their performance and to reduce material costs. The quantity of interest might be calculated from a partial differential equation (PDE) describing the physical phenomena in the machine. \textcolor{red}{This PDE typically depends on some set of parameters, which basically describe geometry and material characteristics of the machine. In industrial workflows they serve as design parameters in numerical models used in the simulation-based engineering process. In the scope of product engineering one often deals with a large amount of conflicting objectives, therefore parametric models are utilized in multi-objective optimizations. They yield solutions on a Parento front~\cite{Baumgartner_2004aa}. In general, the comprehensive assessment of a machine design requires evaluation of product-related physical effects from different domains~\cite{Rosu_2017aa}. However, the treatment of objectives outside the electromagnetic domain is beyond the scope of this paper. High demands on product reliability additionally necessitate the robustness evaluation to be integrated into the optimization process. Since the manufacturing process is always affected by tolerances, these may cause some geometrical and material parameters to become uncertain. As a consequence, the original optimal solution may become suboptimal or even infeasible.} After having chosen a compromise on the Parento front, robust optimization is employed using more sophisticated models. This alleviates the impact of suboptimality by determining solutions that are less influenced by variations. First steps to develop a methodology for robust design were taken by Taguchi (see e.g. \cite{Fowlkes_1995aa}). In this approach, beside the control parameter, the variations are also considered by including noise factors. An overview of different robust optimization methods can be found in \cite{Beyer_2007aa}, where deterministic and evolutionary optimization algorithms are discussed. The former often make use of gradients to find local minima. The advantage is that typically only a few steps are needed to find the minimum \cite{Duan_2013aa}. Robustification is achieved by robust worst-case optimization and the mean-variance approach \cite{Darlington_2000aa}. In the worst-case optimization, the probability density functions (PDFs) of the random parameters are not considered, since the uncertainties are restricted to a bounded uncertainty set \cite{Bertsimas_2011aa}. In the mean-variance approach, the PDFs are considered. The random parameters are described by continuous PDFs. The mean value and the variance are given by integration. The approximation of these integrals is achieved by quadrature. For an example in an optimization setting, one is referred to \cite{Huang_2006aa}. The application of quadrature can be avoided by using approximations based on a linearization of the cost function \cite{Darlington_2000aa}. The introduction of quadrature implies that a lot of \textcolor{red}{evaluations} of the discretized PDE have to be made.

This paper aims at finding the relations between the robust worst-case optimization and the mean-variance approach. The paper addresses equivalences between these methods and proposes efficient solution techniques. To illustrate these findings, the methods are applied to the example of a permanent magnet (PM) synchronous machine (PMSM). 
PMSMs are particularly popular because of their high power density and efficiency. Many aspects of PMSMs have been considered for optimization, e.g., minimization of material costs \cite{Meng_2011aa}. Beside topological rotor shape optimization~\cite{Gangl_2016aa, Putek_2017aa}, also the optimization of the shape of the PMs~\cite{Kim_2005aa, Lukaniszyn_2004aa} has been considered. The PMs are constructed from rare earth elements as is the case for, e.g. NdFeB magnets. The separation of these rare earth elements is environmentally polluting \cite{Binnemans_2013aa}. Therefore in this paper the focus of optimization is on the PMs. In particular, the size of the magnets will be minimized while maintaining a prescribed electromotive force (EMF). \textcolor{black}{In~\cite{Dibarba_2012aa} the mass, including the fraction associated with the PMs, of a PMSM has been reduced while maintaining a desired torque. The optimization was performed using a genetic algorithm. The computational cost has been reduced by relying on 2D finite element machine models and only a posteriori, in the postprocessing step, 3D machine models were calculated. In a subsequent work~\cite{Dibarba_2015aa}, optimization with full 3D models has been conducted, where the authors also used an optimizer with surrogate models to reduce the computational burden.} In~\cite{Li_2017aa} the same goal function is used in the optimization, however, the authors only rely on evolutionary algorithms and do not discuss equivalences between different optimization techniques. To reduce the computational cost they employ kriging~\cite{Lebensztajn_2004aa}. Other ways to optimize electric machines with a reduced order computational method can be found in the literature, e.g.~\cite{Stuikys_2018aa}.

Another approach to reduce the computational costs is Model Order Reduction (MOR). In particular, Proper Orthogonal Decomposition (POD)~\cite{Chatterjee_2000aa} has been \textcolor{black}{shown to be successful. It has been applied for simulating a three-phase transistor~\cite{Henneron_2015ab}, electric machines~\cite{Henneron_2015aa} and high voltage surge arresters~\cite{Schmidthausler_2012aa}. In~\cite{Bontinck_2018aa} an adaptive POD method has been proposed to efficiently model rotating electric machines.}  Using MOR combined with optimization remarkable improvements in speed up have been achieved, e.g. \cite{Gubisch_2017aa,Zahr_2015aa}. \textcolor{black}{Especially when 3D optimizations are performed the use of MOR has shown to increase significantly the efficiency of the optimization procedures e.g.~\cite{Sato_2015aa}.} In the framework of optimizing the magnet of electrical machine\textcolor{red}{s}, POD has been used previously~\cite{Alla_2017aa}. In that paper, the focus is only on worst-case optimization, however the linear approximation for the robust optimization has been extended to a quadratic approximation. The authors opted for MOR since the numerical approximation for robust optimization problems is expensive as it involves solving the PDE numerous times. Another technique for MOR is Reduced Basis (RB)~\cite{Rozza_2008aa,Hess_2013aa}. In \cite{Negri_2013aa, Qian_2017aa} this method has been applied in PDE-constrained optimization.  The same method is used in this work in order to reduce the computational costs.


The paper is structured as follows: In section~\ref{sec:Sens} the forward problem is abstractly stated and different measures for sensitivities are introduced. The various formulations of the optimization problem are discussed in section~\ref{sec:opt_prob}, which are eventually applied to a finite element model in section~\ref{sec:appl}. The results are shown and discussed in section~\ref{sec:Res}. Finally, in the last section conclusions are drawn.

\section{Mathematical Framework}
\label{sec:Sens}
\textcolor{black}{Consider the mathematical model of a machine given by the PDE}
\begin{equation}\label{eq:pde}
	\mathcal{L}\left(u(\pp),\pp\right)=f(\pp),
\end{equation}
on a suitable bounded domain $\mathcal{D}\subset\mathbb{R}^3$ with Dirichlet boundary conditions. The problem depends on the parameters $\pp\in\mathcal{P}$ which we assume to stem from a bounded parametric domain $\parspac\subset\mathbb{R}^{N_\text{P}}$. The solution is given by ${u}$ and $\mathcal{L}$ is an elliptic 
operator with right-hand-side $f$. In the classical \textcolor{black}{2D case of a machine the operator reads} 
\[\mathcal{L}(u(\pp),\pp):=\mathrm{div}\left(a(u,\pp)\nabla u\right).\]
If the dependence of $a$ on $\pp$ is affine, then the solution ${u}$ is analytic in the parameters, cf. \cite{Cohen_2011aa,Babuska_2005ab}. In practice lower regularity may arise \cite{Romer_2016aa} but let us assume that the solution is well-behaved.

The parametrization in terms of $\pp$ is important for optimizing the PDE~\eqref{eq:pde} and for the quantification of uncertainties. Let the linear functional 
\begin{equation}\label{eq:qoi}
	q(\pp):=Q(u(\pp))
\end{equation}
describe a quantity of interest (QoI). It inherits the smoothness of the solution and thus sensitivities can be computed.

\subsection{Local sensitivity analysis}
Let us define, similar to \cite{Cohen_2011aa}, a set of sequences of nonnegative integers
\begin{eqnarray}
\mathcal{F}:= \left\{ \boldsymbol{\alpha}=(\alpha_1,\alpha_2,\ldots): \alpha_i\in\mathbb{N} \wedge \alpha_i \neq 0\right. \nonumber\\
\left.\text{for a finite number of } i \right\},\nonumber
\end{eqnarray}
so that $|\boldsymbol{\alpha}|=\sum_{i\geq1}|\alpha_i|$ is finite if and only if $\boldsymbol{\alpha}\in\mathcal{F}$, then a partial derivative operator $\partial^{\boldsymbol{\alpha}}$ is defined as
\[
\partial^{\boldsymbol{\alpha}}=\frac{\partial^{|\boldsymbol{\alpha}|}}{\partial^{\alpha_1}P_1\ldots\partial^{\alpha_{N_\text{P}}}P_{N_\text{P}}},
\]
with $\boldsymbol{\alpha}\in\mathcal{F}$ supported in $\{1,\ldots,N_\text{P}\}$.
A first tool for sensitivity analysis is relying on sensitivity equations. They require the calculation of the derivatives of the QoI:
\begin{equation}
\label{eq:sensder}
\partial^{\boldsymbol{\alpha}} Q(u(\pp)) = Q({\bf s}^{\boldsymbol{\alpha}}),
\end{equation}
where ${\bf s}^{\boldsymbol{\alpha}} := \partial^{\boldsymbol{\alpha}}u(\pp)$. \textcolor{black}{To increase readability and reduce the number of indices, we will denote, with a slight abuse of notation, the first order sensitivities by} ${\bf s}=[s_1,s_2,\ldots,s_{N_\text{P}}]$ with ${s}_i:=\partial u(\pp)/\partial P_i$.

\textcolor{black}{One defines the reduced space of parameters as
\[
\oparspac := \{\ppexpv\in \parspac | \ppexpv+\ddelta\in\parspac, \ddelta \in\mathcal{U}_\infty \}
\]
to ensure that all perturbations $\ddelta=(\rDelta_1,\ldots,\rDelta_{N_P})^\top$ in the \textcolor{red}{neighborhood} of a point $\ppexpv$ are still in the parameter space $\parspac$ and
\begin{align}
\mathcal{U}_\infty &:= \{\ddelta\in \mathbb{R}^{N_\text{P}}\,|\,\dlow\le \rDelta_i \le \dup,\,i=1,\ldots,N_P\}\nonumber\\ 
&=\{\ddelta\in \mathbb{R}^{N_\text{P}}\,|\,\|D^{-1}\ddelta\|_\infty \le 1\}.
\label{eq:unc_set}
\end{align}
Lower and upper bounds for $\rDelta_i$ are given by $\dlow$ and $\dup$; $D$ is an implicitly defined scaling matrix.} 

\textcolor{black}{Finally, using the above definitions, one can introduce a first order Taylor expansion around $\ppexpv\in\oparspac$
\begin{equation}
\label{eq:tay_exp}
Q(u(\ppexpv+\ddelta))= Q(u(\ppexpv)) + \sum_{i=1}^{N_\text{P}} Q({s}_i)\rDelta_i+\mathcal{O}(\ddelta^2),
\end{equation}
to locally approximate the QoI.}

\subsection{Global sensitivity analysis}
\textcolor{black}{Global sensitivities are defined in a stochastic setting \cite{Sobol_2001aa}}. Therefore it is assumed that $\pp=\pp(\omega)$ are  independent and identically distributed random variables on the probability space $(\Omega,\Sigma,\probspac)$, where $\Omega$ is the set of possible outcomes, $\Sigma$ the sigma algebra and $\probspac$ the probability measure. As a consequence, problem~\eqref{eq:pde} becomes stochastic, i.e., find a stochastic field, $u$, such that it almost surely holds~\cite{Babuska_2005ab}:
\begin{equation}\label{eq:stochpde}
	\mathcal{L}\left(u\left(\pp(\omega)\right),\pp(\omega)\right)=f\left(\pp(\omega)\right).
\end{equation}
\textcolor{black}{By abusing} notation one writes for the unknown $u(\omega)=u\left(\pp(\omega)\right)$ and the QoI $Q(\omega)=Q\left(u(\omega)\right)$, where $\omega$ depicts the stochastic nature of a quantity.

For continuous distributions the expectation value $\EE$ of the QoI or any other stochastic function is defined as 
\[\EE\left[Q\right]=\int_{\Omega}Q(\omega)\;\probspac(\mathrm{d}\omega),\]
and more generally the $k$-th order non-centered moment $\MM$ as
\begin{equation}\label{eq:moments}
\MM^k\left[Q\right]=\int_{\Omega}Q^k(\omega)\;\probspac(\mathrm{d}\omega).
\end{equation}
Global sensitivities can be defined as the first order Sobol coefficients~\cite{Sobol_1990aa,Sobol_2001aa}
\begin{align}
	S_i\left[Q\right]:=\frac{\mathrm{Var}_{i}\left[Q\right]}{\mathrm{Var}\left[Q\right]},
\end{align}
where $\mathrm{Var}[Q]$ is the centered $\MM^2[Q]$ and 
$\mathrm{Var}_i\left[Q\right]:=\mathrm{Var}\left[\EE[Q|P_i]\right]$. 
The meaning of the inner expectation value is that the mean value of $Q$ is taken by considering all $N_P$ parameters as random except for $P_i$, which is kept fixed.

The solution of the probabilistic integrals is typically not available in closed form and one has to carry out numerical quadrature as will be discussed in Section~\ref{sec:appl}.

\subsection{Linearization with respect to the random parameter}
\label{sec:lin_w_mom}
In \cite{Beyer_2007aa} and the references therein, the approaches described in the previous two sections are combined. Assume $\ppexpv=\EE\left[\pp(\omega)\right]$ and $\ddelta'=\ddelta'(\omega)$ to be a random variable, so that $\EE[\rDelta'_i(\omega)]=0$. Then, one finds for the expectation value, $\mu_{Q}=\EE\left[Q(u(\ppexpv+\ddelta'))\right]$, that
\begin{align}
\label{eq:lin_exp}
\mu_{Q}
&= Q(u(\ppexpv)) + \mathcal{O}(\ddelta^{\prime 2}),
\end{align}
where \eqref{eq:tay_exp} is used for the linearization. A similar reasoning can be applied on the variance, $\sigma^2_{Q}=\mathrm{Var}\left[Q(u(\ppexpv+\ddelta'))\right]$,
\begin{align}
\label{eq:lin_var5}
\sigma^2_{Q}&= \mathrm{Var}\left[Q({\bf s})\cdot\ddelta'\right] +\mathcal{O}(\ddelta^{\prime 3})\\&= \sum_{i=1}^{N_\text{P}}Q^2({s}_i)\mathrm{Var}(\ddelta'_i)+\mathcal{O}(\ddelta^{\prime 3}).
\end{align}
Neglecting the higher order terms, the standard deviation (std) becomes
\begin{equation}
\label{eq:lin_var6}
\sigma_{Q}\approx 
\|\mathrm{std}[\ddelta']\circ Q({s}_i)\|_2,
\end{equation}
where $\circ$ depicts the element-wise product. 

\section{Optimization}
\label{sec:opt_prob}
In this section different formulations of the same optimization task are discussed and compared. Two \textcolor{black}{formulations} of the optimization are distinguished: a deterministic and a stochastic one. In the first formulation nominal optimization and classical robust optimization \cite{Hinze_2008aa} are considered using \eqref{eq:pde} as a constraint. The nominal one does not account for deviations on $\pp$ while the latter optimizes the worst case scenario. In the second formulation the optimization of the random PDE is carried out by considering expectation values and standard deviations. Finally similarities between both \textcolor{black}{formulations} are addressed.
\subsection{Deterministic \textcolor{black}{formulation}}
In the deterministic formulation a nominal optimization problem and a robust optimization problem are formulated.
\subsubsection{Nominal optimization}
Let $J_1 : \mathcal{P}\to\mathbb{R}$ depict a cost function which one wants to optimize, so that
\begin{subequations}
\label{eq:opt}
\begin{equation}
 \label{eq:gen_cost_smooth}
 \min_{\ppexpv\in \oparspac} J_1(\ppexpv), \\
\end{equation}
subject to the constraints
\begin{equation}
\label{eq:gen_constraints}
 \textbf{G}_1\left(\ppexpv,q(\ppexpv)\right) \le 0,
\end{equation}
\end{subequations}
where $ \textbf{G}_1$ depicts the collection of functions $G_1^{(j)}$, with $j=1,..., N_G$, with $N_G$ the total number of constraints.
Due to the presence of the QoI in the constraint, the optimization problem is constrained by the PDE \eqref{eq:pde} as discussed in for example \cite{Hinze_2008aa}.

\subsubsection{\textcolor{black}{Robust optimization}}
Let us assume that the parameters $\pp=\ppexpv + \ddelta$ are uncertain within the set $\ddelta\in\mathcal{U}_\infty$, e.g. due to manufacturing imperfections. Hence, a robust counterpart (worst-case) is introduced associated to (\ref{eq:opt}) by considering
\begin{subequations}
\begin{equation}\label{eq:opt_robust}
\min_{\ppexpv\in \oparspac}\max_{\ddelta \in \mathcal{U}_\infty} J_1(\ppexpv+\ddelta),
\end{equation}
subject to
\begin{equation}
\max_{\ddelta \in \mathcal{U}_\infty} \textbf{G}_1(\ppexpv + \ddelta, q(\ppexpv+\ddelta)) \le 0.
\end{equation}
\end{subequations}
This nested optimization problem is hard to solve. Hence, an approximation of the $\max$ problem is utilized. 

\subsubsection{Linearized robust optimization (1-norm)}
By applying a first order Taylor expansion, see e.g. \cite{Diehl_2006aa}, a numerically feasible optimization problem is obtained. Also higher order expansions can be exploited \cite{Lass_2016aa}, however, in this work only linearizations of the cost function and the constraint are considered:
\begin{subequations}
\begin{equation}
\label{eq:RO_lin}
J_1(\ppexpv+\ddelta) \approx  J_1(\ppexpv) + \nabla_{\ppexpv} J_1(\ppexpv)\ddelta 
\end{equation}
and, since there is \textcolor{red}{a} unique solution $u=u(\pp)$ for every admissible point $\pp$, one can introduce the reduced constraints
\[
\tilde{G}_1^{(j)}(\ppexpv):= G_1^{(j)}(\ppexpv+\ddelta, q(\ppexpv+\ddelta)),
\]
which leads to
\begin{align}
\tilde{G}_1^{(j)}(\ppexpv)\approx\;&  G_1^{(j)}(\ppexpv) + \nabla_{\ppexpv} G_1^{(j)}(\ppexpv)\ddelta.
\end{align}
\end{subequations}
Inserting this approximation into the optimization problem, one obtains the linear approximation of the robust counterpart:
\begin{subequations}
\label{eq:opt_robust_1norm}
\begin{equation}\label{eq:opt_robust3}
\min_{\ppexpv\in \oparspac}J_2:= J_1(\ppexpv) + \|D\nabla_{\ppexpv} J_1(\ppexpv)\|_1,
\end{equation}
subject to
\begin{align}
G_2^{(j)}&:=  G_1^{(j)}(\ppexpv, q(\ppexpv)) + \|D\nabla_{\ppexpv} G_1^{(j)}(\ppexpv, q(\ppexpv))\|_1 \le  0,
\end{align}
\end{subequations}
for $j = 1,\ldots,N_G$. The dual norm $||\cdot||_*$ is defined as 
$$
\begin{array}{lrcl}
 \|\cdot\|_*:& \mathbb R^n&\rightarrow&\mathbb R\\
 &\textbf{g}&\mapsto& \|\textbf{g}\|_* := \displaystyle\max_{\textbf{g}\in\mathbb R^n,\|\ddelta\|\le 1} \bf{g}^\top \ddelta.
\end{array}
$$
In this particular case, one can use the property that the dual of $\|D^{-1}\cdot\|_\infty$ is given by $\|D\cdot\|_1$. However, since the norms are not differentiable, this problem is not smooth. To obtain a differentiable problem, slack variables are introduced and one defines a smooth formulation of the optimization problem by
\begin{subequations}
\label{eq:opt_robust_lin}
\begin{equation}
 \label{eq:cost_smooth_robust_lin}
 \min_{\ppexpv\in \oparspac, \boldsymbol{\xi} \in \mathbb R^{N_\text{P}}} J_1(\ppexpv) + \textbf{V}^\top \xi^{(0)}
 \end{equation}
together with the constraints
\begin{equation}
 \label{eq:constraints_robust_lin}
G_2^{(j)}(\ppexpv) + \textbf{V}^\top \xi^{(j)} \le 0
\end{equation}
and
\begin{align}
 \label{eq:constraints_robust_lin_slack}
-\xi^{(0)}& \le D\nabla_{\ppexpv} J_1(\ppexpv)\le \xi^{(0)},\nonumber \\
-\xi^{(j)}& \le D\nabla_{\ppexpv} G^{(j)}_1(\ppexpv)\le \xi^{(j)}, 
\end{align}
\end{subequations}
for all constraints $j=1,\ldots,N_G$, $\boldsymbol{\xi}=(\xi^{(0)}, \ldots,\xi^{(N_\text{P})})$  and $\textbf{V} = (1,\ldots,1)^\top\!\!\in\mathbb R^{N_\text{P}}$. 
This optimization problem can now be efficiently solved numerically. If gradient based methods are used, additionally second-order sensitivities $\mathbf{s}^{\boldsymbol{\alpha}}$, with $|\boldsymbol{\alpha}|\leq 2$ as defined in \eqref{eq:sensder} are required. However, they can be obtained analogously as described previously. Finally, this approach can be generalized to use a quadratic approximation with respect to $\pp$, see~\cite{Lass_2016aa}. 

\subsubsection{Linearized robust optimization (2-norm)}
Albeit less common, one may have chosen the 2-norm instead of the max-norm in the definition of the uncertainty set \eqref{eq:unc_set}. This yields
\begin{equation}
\mathcal{U}_2 :=\{\ddelta\in \mathbb{R}^{N_\text{P}}\,|\,\|D^{-1}\ddelta\|_2 \le 1\}.
\label{eq:unc_set2}
\end{equation}
and changes consequently the norms in \eqref{eq:opt_robust3}. Let us define a reduced parameter space
\[
\oparspactwo := \{\pp\in \parspac | \pp+\ddelta\in\parspac, \ddelta \in\mathcal{U}_2 \}.
\]
The resulting optimization problem reads
\begin{subequations}
\label{eq:opt_robust4}
\begin{equation}
\min_{\ppexpv\in \oparspactwo}J_2:= J_1(\ppexpv) + \|D\nabla_{\ppexpv} J_1(\ppexpv)\|_2,
\end{equation}
subject to
\begin{align} 
G_2&:=G_1^{(j)}(\ppexpv, q(\ppexpv)) + \|D\nabla_{\ppexpv} G_1^{(j)}(\ppexpv, q(\ppexpv))\|_2 \le 0,
\end{align}
\end{subequations}
for $j = 1,\ldots,N_G$. This formulation may not account for optimizing the \textcolor{red}{worst-case} in~\eqref{eq:unc_set} but will still improve its robustness. It does optimize the \textcolor{red}{worst-case} for \eqref{eq:unc_set2} and has been used in \cite{Diehl_2006aa}.
\subsection{Stochastic \textcolor{black}{formulations}}
In the stochastic formulation the cost function and the constraints \eqref{eq:stochpde} are stochastic due to the uncertainties on $\pp$, which are given by random distributions. Nonetheless, the problem can be formulated in terms of the stochastic moments.
\subsubsection{Nominal optimization}
Recalling that $\pp=\ppexpv+\ddelta^{\prime}$, the cost function and its constraints can now be defined in terms of expectation values, as in \cite{Sundaresan_1995aa},
\begin{subequations}
\label{eq:gen_opt_uq}
\begin{equation}
 \label{eq:gen_cost_uq}
 \min_{\ppexpv\in \oparspactwo} J_3(\pp):=\EE\left[J_1(\pp)\right],\\
\end{equation}
subject to
\begin{equation}
\label{eq:gen_const_uq}
 G^{(j)}_{3}\left(\pp,q(\pp)\right):= \EE\left[G_{1}^{(j)}(\pp,q(\pp))\right]\le 0.
\end{equation}
\end{subequations}
This optimization problem is again deterministic and can be solved by using the same techniques as for the nominal deterministic optimization, however, with increased computational costs for approximating the probabilistic integrals.

\subsubsection{Robust optimization}
The approach for robust optimization with stochastic quantities is similar to the description mentioned in \cite{Darlington_2000aa}, where in the cost function the expectation value and the variance are considered. They are weighted by a, so-called, risk aversion parameter. In this work this parameter has been normalized. In \cite{Huang_2006aa} also the constraints are robustified by including the standard deviations. Combining the two ideas, results in
\begin{subequations}
\label{eq:gen_opt_uq_rob}
\begin{align}
 \label{eq:gen_cost_uq_rob}
 \min_{\ppexpv\in \oparspactwo}J_4(\pp):=\;&\EE\left[J_1(\pp)\right]+ \lambda\hspace{0.5em}\std\left[J_1(\pp)\right],
\intertext{subject to}
\label{eq:gen_const_uq_rob}
 G_{4}^{(j)}(\pp,q(\pp)) :=\;& \EE\left[G_{1}^{(j)}(\pp,q(\pp))\right]\\
 &+ \lambda\hspace{0.1em}\std\left[G_{1}^{(j)}(\pp,q(\pp))\right] \le 0,\nonumber
\end{align}
\end{subequations}
where $\lambda>0$ can be interpreted as a weighting factor, similar to $D$ in \eqref{eq:unc_set}.

If one wants to increase the stochastic accuracy, one can easily take higher order moments into account, see e.g. \cite{Tiesler_2012aa}.
\subsubsection{\textcolor{black}{Linearized robust optimization}}
Applying the linearization of \eqref{eq:lin_exp} and \eqref{eq:lin_var6} to the cost function and its constraints, leads to
\begin{subequations}
\label{eq:gen_opt_uq_rob_lin}
\begin{align}
 \label{eq:gen_cost_uq_rob_lin}
 J_5(\ppexpv)=\;&\EE\left[J_1(\ppexpv+\ddelta')\right]+\lambda\;\mathrm{std}[J_1(\ppexpv+\ddelta')]\nonumber\\ 
 \approx\;& J_1(\ppexpv)+\lambda\|\mathrm{std}[\ddelta']\circ \nabla_{\ppexpv}J_1(\ppexpv)\|_2.
\end{align}
and
\begin{align}
\label{eq:gen_const_uq_rob_lin}
 G_{5}^{(j)}(\ppexpv,q(\ppexpv))=\;&\EE\left[G_{1}^{(j)}(\ppexpv+\ddelta',q(\ppexpv+\ddelta'))\right]\nonumber\\
 &+ \lambda\;\mathrm{std}\left[G_{1}^{(j)}(\ppexpv+\ddelta',q(\ppexpv+\ddelta'))\right]\nonumber\\
 \approx\;& G_{1}^{(j)}(\ppexpv,q(\ppexpv))\nonumber\\
 &+\lambda\|\mathrm{std}[\ddelta']\circ \nabla_{\ppexpv}G^{(j)}_1(\ppexpv,q(\ppexpv))\|_2.
\end{align}
\end{subequations}
This reasoning can be extended to a quadratic approximation~\cite{Alexandrian_2017aa}. Finally, considering the linearization above unveils the following equivalence
\begin{theorem} Let the assumptions
	(i) the QoI $Q$ from \eqref{eq:qoi} is linear in $u$, such that $\partial Q/\partial u=$const,
	(ii) perturbations are independent and identically distributed and $\rDelta'_i$ is symmetric around 0, 
	hold true, then, by choosing $\lambda=\frac{D_{ii}}{\mathrm{std}[\rDelta_i']}$, the linearized UQ optimization (i.e. \eqref{eq:gen_opt_uq_rob_lin}) is equivalent to the worst-case optimization using the 2-norm (i.e. \eqref{eq:opt_robust4}).
\end{theorem}

In contrast to related results in the literature this theorem focuses on the case of PDE-constraint optimization, discusses the assumptions rigorously and translates all parameters among the two approaches. 

For \eqref{eq:gen_opt_uq_rob_lin} first order quadrature is sufficient for determining the standard deviations. Doing this quadrature for~\eqref{eq:gen_opt_uq_rob} leads to a different linear approximation in parameter space. However, both approaches would be affected by the curse of dimensionality ($2^{N_{\pp}}$ quadrature points). Using \eqref{eq:opt_robust4} offers a \textcolor{black}{computationally} cheap alternative, especially when the (higher order) derivatives are available. On the other hand the embedding of \eqref{eq:opt_robust4} into the stochastic framework allows the development of hybrid optimization algorithms which adaptively switch among methods depending on the accuracy required for the probabilistic integrals. 

\subsection{\textcolor{black}{Algorithms for optimization}}
\textcolor{black}{If sensitivities \eqref{eq:sensder} are available, the optimization problem can be straightforwardly solved by gradient-based methods, e.g. Sequential quadratic programming (SQP) with damped Broyden-Fletcher-Goldfarb-Shanno (BFGS) updates \cite{Nocedal_2006aa}. If no derivatives are available one can still use deterministic algorithms, for example bound optimization by quadratic approximation (BOBYQA)~\cite{Powell_2009aa} or genetic, swarm and evolution based algorithms, e.g. \cite{Li_2017aa, Sizov_2013aa}.}

\textcolor{black}{The idea of particle swarm optimization (PSO) is that a set of $\Npar$ particles are moving through the admissible design space spanned by the parameters~\cite{Kennedy_1995aa}. At the position of every particle the objective function is evaluated. In the next iteration step ($k+1$) the position of every particle is updated taking into the previous best values in the history of every particle and the best value of the swarm. The best sets are denoted by $\phat_n$ for particle $n$ and $\phat_{\mathrm{sw}}$ for the swarm. The ``velocities''of every member of the swarm are updated by}
\textcolor{black}{
\begin{align}
\nonumber
\mathbf{v}_{k+1,n}=
\omega_0 \mathbf{v}_{k,n} &+\omega_1\mathbf{R}_1\left(\phat_n-\pp_{k,n}\right)\\
\label{eq:PSO}
&+\omega_2\mathbf{R}_2\left(\phat_{\mathrm{sw}}-\pp_{k,n}\right),
\end{align}
with $\omega_i$ ($i\in\{1,2,3\}$) swarm characteristic constants and $\mathbf{R}_j$ ($j\in\{1,2\}$) two random diagonal matrices where the diagonal entries for both matrices $r_{kk}$ are independently and uniformly chosen from the interval $[0,1]$. These matrices mimic the free will of the swarm. The first term in~\eqref{eq:PSO} tries to maintain a part of the current velocity, the second term convinces the particle to head to the best found point $\phat_n$ and the last term tries to drive the particle into the direction of the swarm's best found point. When a particle is leaving the admissible set it is projected back on the boundary of the set. As initialization the particles are randomly and uniformly distributed over the set and $\mathbf{v}_{0,n}$ is set to zero for all particles. There are three stopping criteria for the algorithm, namely:
\begin{enumerate}
\item the maximum number of iterations has been reached;
\item the majority of the particles are close enough to $\phat_{\mathrm{sw}}$, so that
\[
\frac{1}{\Npar}\sum_{n=1}^{\Npar} \|\phat_{\mathrm{sw}}-\pp{k,n}\|_2<\epsilon
\]
is fulfilled with $\epsilon$ a \textcolor{red}{user-defined} accuracy;
\item there is no further improvement in $\phat_{\mathrm{sw}}$ over $N_{\mathrm{stall}}$ consecutive iterations.
\end{enumerate} }


\section{Application}
\label{sec:appl}
The above mentioned optimization procedures are applied to a 3-phase 6-pole PMSM with buried magnets. 
The goal of the optimization is to reduce the amount of PM material while maintaining a QoI, here, in particular, the electromotive force.

\subsection{Machine description}
\begin{figure}
      \centering
      \includegraphics[width=0.45\textwidth]{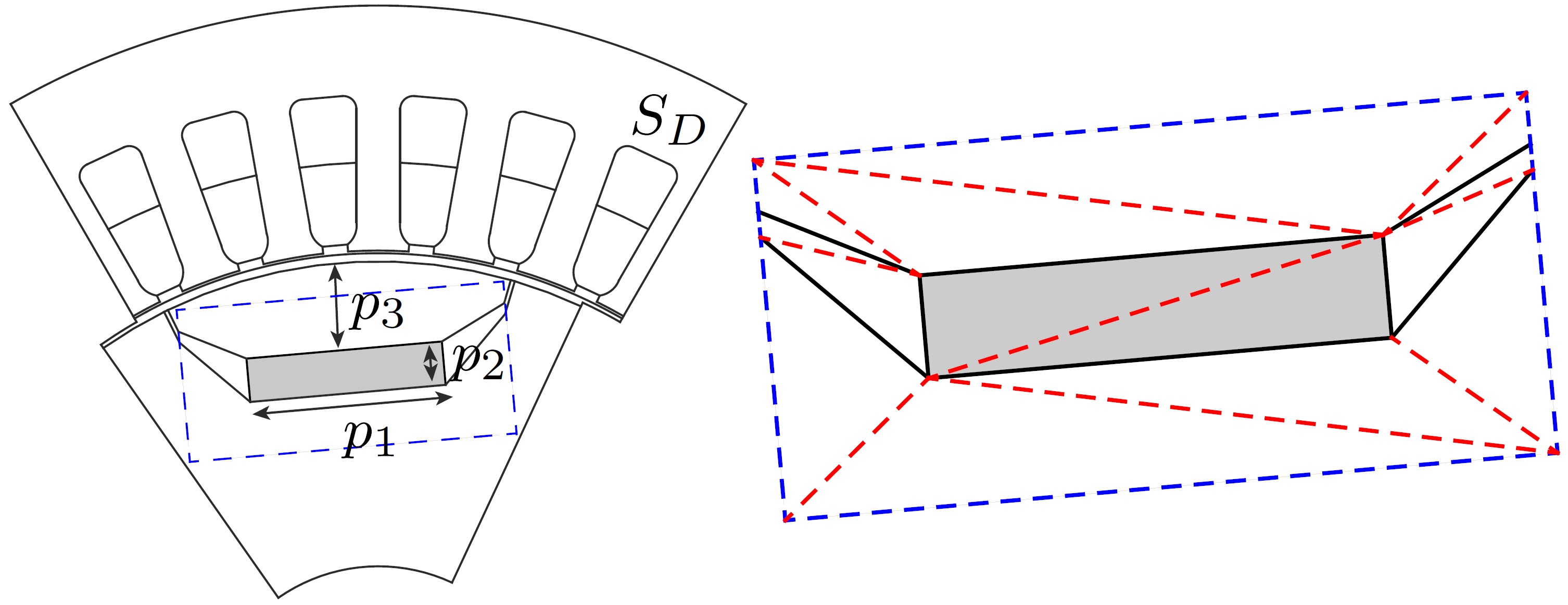}
\caption{The cross-section of one pole of the machine with the magnet depicted in gray and the region of the affine decomposition indicated by the dashed box. On the right hand side, the triangulation into $N_L $ subdomains is shown by the dashed lines.}
\label{fig:aff_dec}
\end{figure}

One pole of the 3-phase PMSM is shown in Fig.~\ref{fig:aff_dec}.  The winding is a double layered winding with two slots per pole per phase. The laminated steel of the machine is modeled with zero conductivity and a relative permeability of $\mu_r = 500$. The machine is based on the model described in \cite{Pahner_1998aa}. The magnets have a width~$\pone$ and height $\ptwo$. They are buried in the rotor at a depth~$\pthree$.  The machine is subjected to an optimization procedure in which the amount of PM material is reduced while maintaining a prescribed EMF $E_\mathrm{d}$. 

\textcolor{black}{For cylindrical machines with lengths comparable to or larger than their diameters, 2D field models are typically used to calculate the major machine parameters up to an accuracy that is sufficient in practice. Occurring 3D effects (such as, e.g., the resistance and leakage inductance of the end-winding parts and the different lengths of stator and rotor lamination stacks) are included as additional lumped elements at the circuit level or by adaptation of the material parameters. In earlier design steps, such adapted 2D models are preferred over 3D models because of the according smaller computation times. Commonly afterwards 3D high-fidelity simulations are carried out, e.g.~\cite{Dibarba_2012aa,Pels_2015aa}.}

\subsection{Finite element model}
To calculate the EMF of every configuration, one has to solve Maxwell's equations. The magneto\-static formulation is sufficient to model this type of electric machines~\textcolor{black}{\cite{Nasar_1993aa}}. This implies that eddy currents and displacement currents are neglected. To calculate the magnetic vector potential (MVP) $\MVP(x,y,z,\pp)$, one has to solve the PDE
\begin{equation}\label{eq:mqscont}
    \vec{\nabla}\times\left(\nu(\pp)\vec{\nabla}\times\MVP(\pp)\right) =\vec{J}_{\rm src}-\vec{\nabla}\times\vec{H}_{\rm pm}(\pp),
\end{equation}
with adequate boundary conditions. $\nu(\pp)=\nu(x,y,z,\pp)$ depicts the reluctivity, $\vec{J}_{\rm src}(x,y,z)$ is the source current density and $\vec{H}_{\rm pm}(x,y,z,\pp)$ is the PM's source magnetic field strength. The magnetization current density induced by the PMs, i.e., $\vec{\nabla}\times\vec{H}_{\rm pm}(\pp)$. 
Nonlinear saturation curves are not discussed, but could be considered easily. In the 2D planar case, the ansatz
\[
\MVP(\pp)\approx\sum_{j=1}^{N_D} u_j(\pp) \vec{w}_j(x,y)=\sum_{j=1}^{N_D} u_j(\pp) \frac{N_j(x,y)}{l_z} \vec{e}_z
\]
is made, where $N_D$ is the total number of degrees of freedom, $\vec{w}_j(x,y)$ are the edge shape functions related to the nodal finite elements $N_j(x,y)$ and $l_z$ is the length of the machine.
The Galerkin procedure leads to the system of equations 
\begin{equation}
\label{eq:sys_eq}
\mathbf{K}_\nu (\pp) \mathbf{u}(\pp) =\mathbf{j}_{\rm src}+\mathbf{j}_{\rm pm}(\pp)\end{equation} 
with 
\begin{align}
K_{\nu(\pp),i,j}(\pp)&=\int_{V_D} \nu(\pp) \vec{\nabla}\times\vec{w}_i\cdot\vec{\nabla}\times\vec{w}_j\; \text{d}V,
\\
j_\text{src,i}&=\int_{V_D}\vec{J}_{\rm src}\cdot \vec{w}_i\;\text{d}V,
\\
j_\text{pm,i}(\pp)&= -\int_{V_D}\vec{H}_{\rm pm}(\pp) \cdot \vec{\nabla}\times\vec{w}_i\;\text{d}V,
\end{align}
where $V_D=S_D \times[0,l_z]$ is the computational domain and $S_D$ the machine's cross-section \cite{ Monk:2003aa}.

Solving \eqref{eq:sys_eq} gives the MVP from which the EMF $E_0$ can be calculated by using the loading method~\cite{Rahman_1991aa}. It corresponds essentially to an FFT. This means that the EMF can be computed during post-processing. Using the notation above, our \textcolor{red}{QoI} reads $$E_0(\mathbf{u}(\pp))=q(\pp).$$

\subsection{Affine decomposition}
\label{sec:aff_dec}
Since during the optimization different configurations of the magnet position in the machine will be considered, one wants to resort to a computationally fast model and avoid remeshing in order to reduce numerical noise. Therefore an affine decomposition (see e.g. \cite{Rozza_2008aa, Alla_2017aa}) is introduced into the model. Hence, a region around the permanent magnet (Fig.~\ref{fig:aff_dec}) is subdivided into $N_L$ triangular subdomains. The finite system matrix from~\eqref{eq:sys_eq} can be rewritten as
\begin{equation}
\label{eq:sys_eq_aff}
\mathbf{K}_{\nu}(\pp) =\mathbf{K}_{\nu}^{0}+\sum_{\ell=1}^{N_L} \vartheta^{\ell}(\pp)\mathbf{K}_{\nu}^{\ell},
\end{equation}
where $\mathbf{K}_{\nu}^{0}$ represents the system matrix  for the domain outside the dashed box and $\mathbf{K}_{\nu}^{\ell}$ represents the matrices corresponding to the subdomains. The weights $\vartheta^{\ell}$ inherit the dependency on $\pp$ and can be computed analytically. The notation in \eqref{eq:sys_eq_aff} is compact for
\begin{eqnarray}
\label{eq:comp_not}
\vartheta^{\ell}(\pp)\mathbf{K}_{\nu}^{\ell} &:=& \vartheta^{\ell}_1(\pp) \mathbf{K}^{\ell}_{xx}+\vartheta^{\ell}_2(\pp)\mathbf{K}^{\ell}_{yy}\nonumber\\
& &+\vartheta^{\ell}_3(\pp) \mathbf{K}^{\ell}_{xy}+\vartheta^{\ell}_4(\pp) \mathbf{K}^{\ell}_{yx},\nonumber
\end{eqnarray}
where $\mathbf{K}_{\nu}^{\ell}$ has been decomposed into four submatrices $\mathbf{K}^{\ell}_{xx}$, $\mathbf{K}^{\ell}_{yy}$, $\mathbf{K}^{\ell}_{xy}$ and $\mathbf{K}^{\ell}_{yx}$. The subindices indicate the partial derivatives of the nodal shape functions. 
An additional advantage of this method is that these matrices can be precomputed and matrix assembly can be avoided. An analogue decomposition is made for the right-hand-side of~\eqref{eq:sys_eq}. \textcolor{black}{The application of the affine decomposition can be extended to 3D~\cite{Hess_2013aa}.}

\subsection{Stochastic quadrature}
In the stochastic setting, the probabilistic dimension of the PDE~\eqref{eq:stochpde} is sampled by using collocation, \cite{Xiu_2010aa}. The integrals \eqref{eq:moments} are approximated by

\begin{equation}
	\MM^k\left[q\right]\approx \sum_{n=1}^N w_n q^k(\pp_n),\nonumber
\end{equation}
where $N$ are the number of samples and the weights $w_n$ and evaluation points $\pp_n$ are method specific. In this paper stochastic quadrature (SQ) and Monte Carlo (MC) are used. In the former case the weights and sample points $\pp_{n}$ are chosen according to the quadrature rules. For the latter, all samples have the same weight $1/N$ but the sample points itself are chosen randomly. The main advantage of SQ  over MC is the fast convergence for low dimensional problems \cite{Xiu_2010aa}. 

The computational drawback of all quadrature-based approaches is the need of many evaluations of the PDE \eqref{eq:pde}. This is computationally expensive and therefore MOR is desirable. 
\subsection{Model order reduction}
\label{sec:mor}
To speed up the computations the reduced basis method is used \cite{Rozza_2008aa}. The idea is to project the high dimensional problem to a lower dimensional space. One looks for a solution in the reduced space $\{\psi_1,\ldots,\psi_d\}$ of rank $d$. This space is constructed during an \textit{offline} phase. The basis $\psi_i$ of the space is computed by a greedy algorithm. With the help of an error estimator, it is decided for which $\pp$ in a training set, a solution is computed and added as $\psi_i$ in order to obtain a maximal error reduction. For every $\psi_i$ the high dimensional problem is solved and then orthonormalized with respect to the current subspace by using the Gram-Schmidt process. The algorithm is stopped when a predefined accuracy is reached.

Finally one obtains a reduced basis of rank $d$ that is sufficient to capture the dynamics in the parameter space. With the Galerkin ansatz one retrieves 
\[{\bf u}^d(\pp):=\sum_{i=1}^d \tilde{\mathbf{u}}_i(\pp)\psi_i = \Psi \tilde{\mathbf{u}}(\pp).\]
Multiplying equation \eqref{eq:sys_eq} from the left with $\Psi^{\top}$, utilizing the affine decomposition introduced in \eqref{eq:sys_eq_aff} and inserting ${\bf u}^d(\pp)$, one retrieves the reduced order model
\begin{align*}
\Psi^{\top}\left(\mathbf{K}^{0}+\sum_{\ell=1}^{L}\vartheta^\ell(\pp) \mathbf{K}^\ell\right)\Psi\tilde{\mathbf{u}}(\pp)\\= \Psi^{\top}\left(\mathbf{j}^{\rm out}+\sum_{\ell=1}^{L}\vartheta^\ell(\pp) \mathbf{j}^\ell\right).
\end{align*}
Since the problem is linear one obtains 
\begin{align*}
\left(\underbrace{\Psi^{\top}\mathbf{K}^{0}\Psi}_{\tilde{\mathbf{K}}^{\rm out}}+\sum_{\ell=1}^L\vartheta^\ell(\pp)\underbrace{\Psi^{\top}\mathbf{K}^\ell\Psi}_{\tilde{\mathbf{K}}^\ell}\right)\tilde{\mathbf{u}}(\pp)\\=\underbrace{\Psi^{\top}\mathbf{j}^{\rm out}}_{\tilde{\mathbf{j}}^{\rm out}}+\sum_{\ell=1}^L\vartheta^\ell(\pp)\underbrace{\Psi^{\top}\mathbf{j}^\ell}_{\tilde{\mathbf{j}}^\ell}.
\end{align*}
Note that all quantities with tildes are of dimension $d \ll N_D$. They can also be precomputed, except for $\tilde{\mathbf{u}}(\pp)$. The reduced system of equations can now be solved very efficiently during the \textit{online} phase. The system can also be set up for different values of $\pp$ without the need for high dimensional operations. This benefit is obtained thanks to the particular affine decomposition introduced in Section~\ref{sec:aff_dec}.

The error on the solution can be bounded by a posteriori error estimates (see for example \cite{Rozza_2008aa}). Define a norm on a reference geometry by  $\|\mathbf{v}\|^2_{\pp_\mathrm{ref}} = \mathbf{v}^\top \mathbf{K}(\pp_\mathrm{ref}) \mathbf{v}$. Then the error estimator is given by 
\[\|\mathbf{u}(\pp)-\mathbf{u}^d(\pp)\|_{\pp_\mathrm{ref}}\leq \Delta_\mathbf{u}(\pp):= \frac{\| \mathbf{r}(\pp)\|_{\pp_\mathrm{ref}*}}{\alpha(\pp)}. \] 
The residual $\mathbf{r}$ is defined by $\mathbf{r}=\mathbf{K}(\pp)\mathbf{u}^d-\left(\mathbf{j}_{\rm src}+\mathbf{j}_{\rm pm}\right)$. The dual norm is depicted by $\|\mathbf{v}\|^2_{\pp_\mathrm{ref}*} = \mathbf{v}^\top \mathbf{K}(\pp_\mathrm{ref})^{-1} \mathbf{v}$. Due to the affine decomposition and our choice for the norm, the coercivity constant $\alpha(\pp)$ of $\mathbf{K}(\pp)$ can be computed by the "min $\Theta$" approach
\[\alpha(\pp)=\min_{ \ell\in\{1,\ldots, N_{\ell}\}}\frac{\vartheta^{\ell}(\pp)}{\vartheta^{\ell}(\pp_\mathrm{ref})}.\] 
Furthermore, an error estimator for the sensitivity can be derived with similar ingredients. Detailed information can be found in ~\cite{Oliveira_2007ab}. The error estimator can be decomposed in the offline/online framework. Hence, in the online phase the evaluation of the error estimator does not rely on high dimensional operations.

Problems with a high sensitivity with respect to the parameter pose major challenges when the reduced basis method is applied. Especially problems involving geometry transformations can be very sensitive to the parameters and lead to large reduced order models since different
phenomena have to be captured. To further reduce the computational cost during the online phase and to keep the dimension of the reduced order models small, a 'Dictionary' of models is generated. This is obtained by dividing the parameter space into $N_\mathrm{Par}=40$ partitions. Each partition represents a cube $Q_{ijk}$ in the parameter space, which is defined as $Q_{ijk}=[t_i^{(1)}, t_{i+1}^{(1)}]\times[t_j^{(2)}, t_{j+1}^{(2)}]\times[t_k^{(3)}, t_{k+1}^{(3)}]$ with $t_i^{(l)}\in p^{(l)}$, with $p^{(1)}=[0.5, 3, 6, 8, 11, 13, 16, 18, 21, 23, 26.5]$, $p^{(2)}=[0.5, 7.5, 10.5]$ and  $p^{(3)}=[4.5, 7.5, 14.5]$. For each partition a separate reduced order model is generated. In the online phase for a given parameter $P$ the associated partition is determined and the corresponding reduced model utilized. This approach allows us to obtain low dimensional models that can be evaluated rapidly. A similar approach has been investigated in \cite{Haasdonk_2011aa} where a strategy using adaptive partitioning was developed. Optionally, in the presented approach the offline phase can be accelerated significantly by using parallel computing, since the partitioning is chosen fixed and the reduced order models in the different partitions can be computed independently.

\subsection{Optimization procedure}
The modeling of the machine is carried out in 2D, the optimization considers the reduction of the surface $S=\pone\ptwo$ instead of the volumes of the magnet. The depth of the magnet $\pthree$ is chosen to be a free parameter which is also changed during the optimization process. 

\subsubsection{Problem definition}
In the deterministic nominal optimization the cost function \eqref{eq:gen_cost_smooth} is given as
\begin{subequations}
\label{eq:opt1}
\begin{equation}
 \label{eq:cost_smooth}
 \min_{\ppexpv\in \mathbb R^3} J_1(\ppexpv) := \pexpv_1 \pexpv_2\\
\end{equation}
subject to
\begin{equation}
\label{eq:constraints}
 G_1(\ppexpv,q(\ppexpv)) := \left(\begin{array}{c}
   \pplow_1 - \pexpv_1\\
   \pplow_2  - \pexpv_2\\
   \pplow_3  - \pexpv_3\\
   \pexpv_3 - \ppup_3\\
  \pexpv_2 + \pexpv_3 - 15\\
  3\pexpv_1 - 2\pexpv_3 - 50\\  
  E_\mathrm{d} - E_0({\bf u}(\ppexpv))
 \end{array}\right) \le 0.
\end{equation}
\end{subequations}
The first four constraints are related to the lower and upper bounds of $\pp$: $(\pplow_1,\pplow_2,\pplow_3) = (1,1,5)$ and $(\ppup_1, \ppup_2, \ppup_3) =(\infty, \infty, 14)$. The fifth constraint ensures the validity of the affine decomposition (no intersections). Only a sub-domain of the geometry is considered. Hence, it is required to stay in that region. The sixth constraint is a design constraint, enforcing that each PM has to have a certain distance to the rotor's surface, meaning that the depth of the magnet is linked to its width. The last constraint is the requirement to fulfill the prescribed EMF and since it is calculated from \eqref{eq:sys_eq}, the optimization problem actually has a PDE constraint. For \eqref{eq:opt_robust_lin}~and \eqref{eq:opt_robust4} the uncertainty set is chosen to be $D = \mbox{diag}((\dup - \dlow)/2)$, where in our numerical experiments $-\dlow = \dup=\dbound$ and the value of $\dbound$ is increased from 0 to $\unit{0.2}{mm}$. 

In the stochastic \textcolor{black}{formulation}, i.e.~\eqref{eq:gen_opt_uq} and \eqref{eq:gen_opt_uq_rob}, it is assumed that the components of $\pp$ are independently uniformly distributed: 
\begin{equation}
\label{eq:distr}
\pp\sim\mathcal{U}(\ppexpv+\ddlow,\ppexpv+\ddup),
\end{equation}
where $\ppexpv=\mathbb{E}\left[\pp\right]$. Since $\pone$ and $\ptwo$ are independent random variables, \eqref{eq:gen_cost_uq} can be written as
\begin{subequations}
\label{eq:opt_uq}
\begin{equation}
 \label{eq:cost_uq}
 \min_{\pp\in \mathbb R^3} J_3= \EE\left[\pone\right] \EE\left[\ptwo\right],\\
\end{equation}
subject to
\begin{equation}
\label{eq:const_uq}
 G_3 := \left(\begin{array}{c}
   \pplow_1 - \EE\left[\pone\right]\\
   \pplow_2 - \EE\left[\ptwo\right]\\
   \pplow_3 - \EE\left[\pthree\right]\\
   \EE\left[\pthree\right] - \ppup_3\\
  \EE\left[\ptwo\right] +\EE \left[\pthree\right] - 15\\
  3\EE\left[\pone\right] - 2\EE\left[\pthree\right] - 50\\  
  E_\mathrm{d} - \EE\left[E_0({\bf u}(\ppexpv))\right]
 \end{array}\right) \le 0,
\end{equation}
\end{subequations}
where the notation has been shortened so that $J_i=J_i(\pp)$ and $G_i= G_i(\pp,q)$. For the robust counterpart in the stochastic setting, the same convention for notation is applied:
\begin{subequations}
\label{eq:opt_uq_rob}
\begin{equation}
 \label{eq:cost_uq_rob}
 \min_{\pp\in \mathbb R^3}J_4:=\EE\left[\pone\right] \EE\left[\ptwo\right]+ \lambda\hspace{0.5em} \std\left[\pone\ptwo\right],\\
\end{equation}
subject to
\begin{align}
\label{eq:const_uq_rob}
 G_4 := &\left(\begin{array}{c}
   \pplow_1 - \EE\left[\pone\right]+\lambda\std\left[\pone\right]\\
   \pplow_2 - \EE\left[\ptwo\right]+\lambda\std\left[\ptwo\right]\\
   \pplow_3 - \EE\left[\pthree\right]+\lambda\std\left[\pthree\right]\\
   \EE\left[\pthree\right]-\lambda\std\left[\pthree\right] - \ppup_3\\
  \EE\left[\ptwo\right] +\EE \left[\pthree\right] -\lambda\std\left[\cdot\right] - 15\\
  \EE\left[3\pone\right] - \EE\left[2\pthree\right]+\lambda\std\left[\cdot\right] - 50\\  
  E_\mathrm{d} - \EE\left[E_0\right]+\lambda\std\left[E_0({\bf u}(\ppexpv))\right]
 \end{array}\right)
 \le 0. 
\end{align}
\end{subequations}
The sampling in both cases is carried out using SQ and MC. For SQ a tensor grid of $5\times5\times5$ is constructed and a Gau\ss-Legendre quadrature is applied since only uniform distributions are considered \cite{Xiu_2010aa}. For MC $N_\text{MC}=5000$ random samples are generated such that the estimated Monte Carlo error is below $1\%$ for all optimizations.

\subsubsection{Optimization algorithms}
\textcolor{black}{Since first and second order sensitivities are available due to the affine decompositions, SQP with BFGS is used \cite{Nocedal_2006aa,Lass_2016aa}.} \textcolor{black}{The algorithm stops when an \textcolor{red}{accuracy} of $10^{-3}$ is obtained or the maximum number of 10 iterations is reached.} \textcolor{black}{On the other hand, stochastic algorithms are often applied in machine optimization. Therefore the deterministic algorithms discussed above are compared to the genetic algorithm (GA) implemented in MATLAB\textsuperscript{\textregistered}  and to an in-house implementation of PSO. The robustification is applied by considering $J_4$ as a cost function and with the constraints $G_4$. To keep the computational cost reasonable, only SQ, with a $5\times5\times5$-grid was used to determine the expectation values and the standard deviations. The computations are parallelized over 4 processor cores.}

\textcolor{black}{For the PSO, the number of particles at every iteration is set to 50, $N_\mathrm{stall}=15$ and the maximum number of iterations is put to 100. The tolerance $\epsilon$ has been put to $10^{-6}$. The constants in \eqref{eq:PSO} are set to $\omega_0=0.5$ and $\omega_1=\omega_2=1.49$. }

\section{Results and Discussion}
\label{sec:Res}
In Table~\ref{tab:res_nomopt} the results for the nominal optimization \textcolor{red}{are} shown. The results of the different robust optimization procedures are depicted in Table~\ref{tab:res_robopt}, where a maximal deviation $\dbound=\unit{0.2}{mm}$ was considered. All calculations were run on a $\unit{16}{GB}$ RAM Intel\textsuperscript{\textregistered} Core\textsuperscript{\texttrademark} with i7-3820K processors ($\unit{3.60}{GHz}$).

\begin{figure}%
\centering
\subfloat[][Initial]{\includegraphics[width=0.25\textwidth]{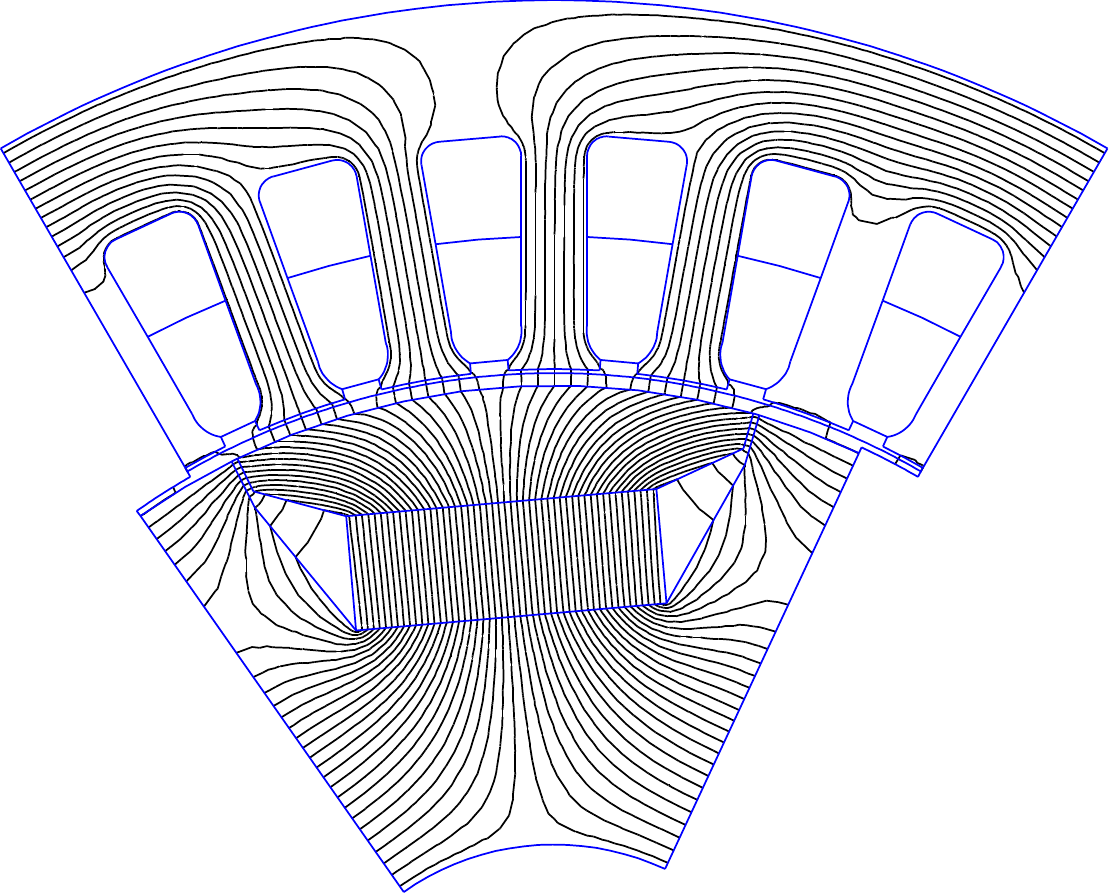}}%
\subfloat[][\textcolor{black}{Nominal (i)}]{\includegraphics[width=0.25\textwidth]{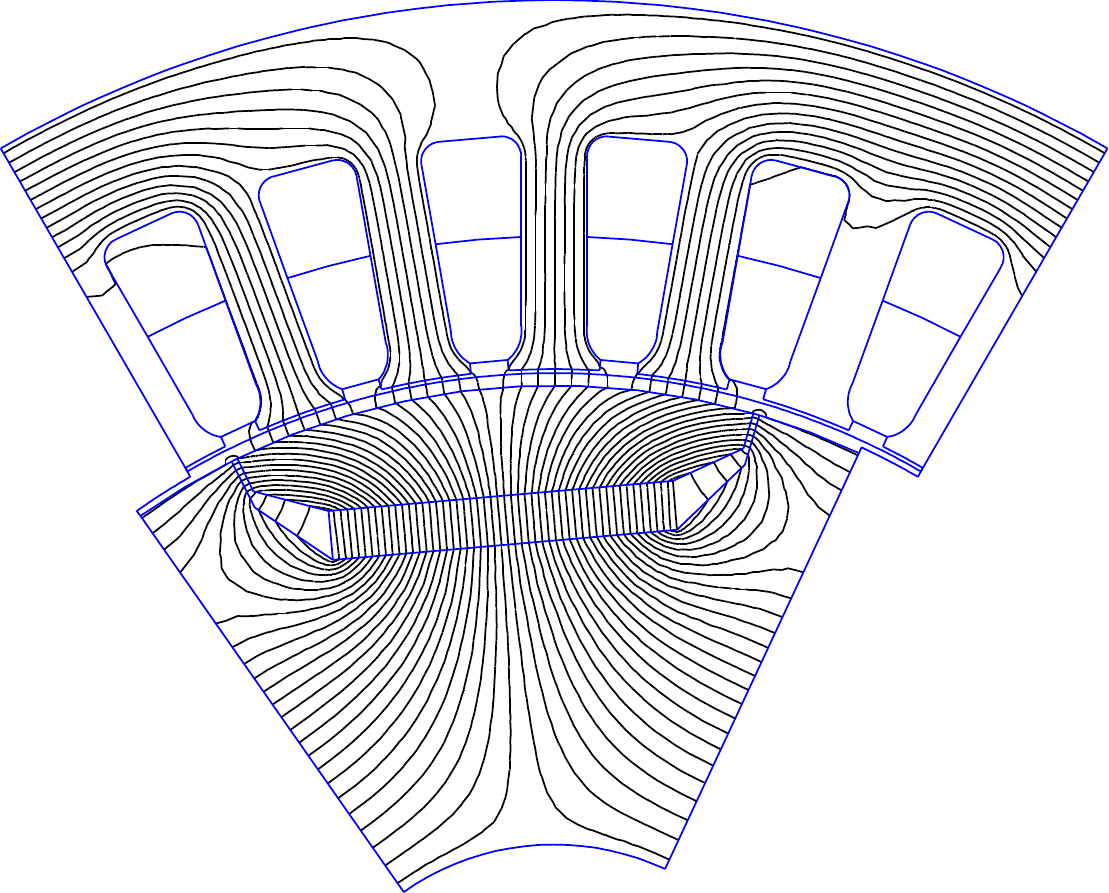}}\newline
\subfloat[][\textcolor{black}{Linearized (i)}]{\includegraphics*[width=0.25\textwidth]{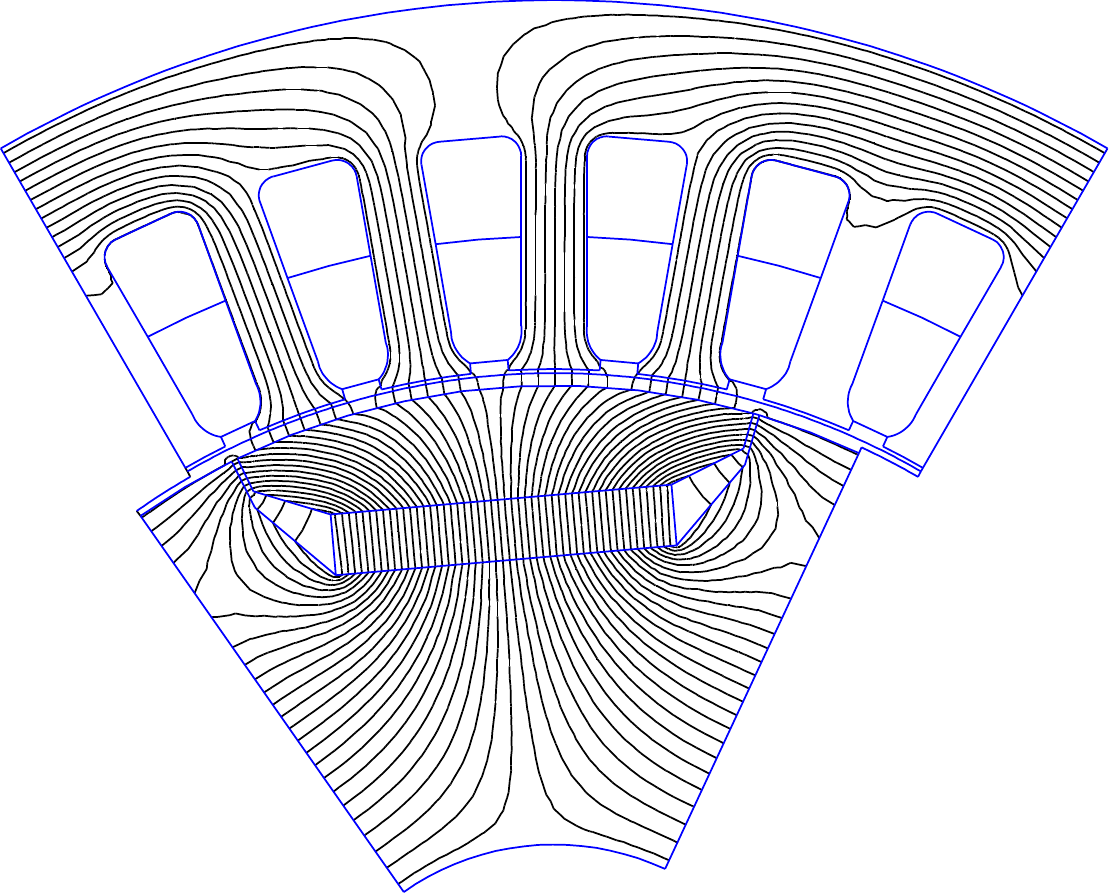}}%
\subfloat[][\textcolor{black}{Robust (v)}]{\includegraphics[width=0.25\textwidth]{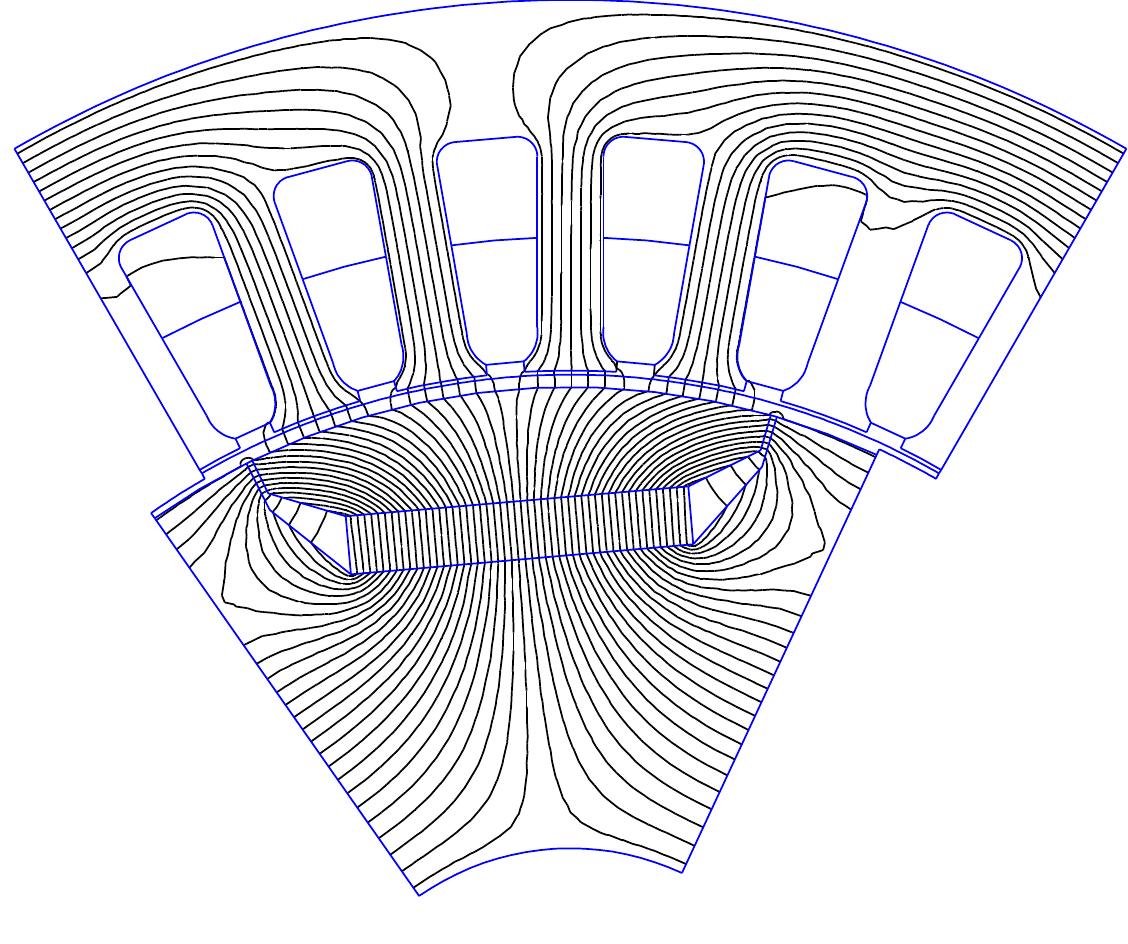}}
\caption[]{Optimized PMSM design according to three different algorithms.}%
\label{fig:results}%
\end{figure}

A visualization of the machines is depicted in Fig.~\ref{fig:results}.  One sees that the size of the PM has decreased substantially while maintaining the desired EMF. For the settings with nominal optimization the three methods result in comparable PM sizes. The volume has been reduced by more than $\unit{50}{\%}$, which indicates a very good improvement (and that the initial guess was poor). Since the offline phase for the construction of the reduced basis takes $\unit{234}{s}$, there is no MOR applied for Nominal (i), Linearized (i) and Linearized (ii). The UQ deterministic settings on the other hand require many more evaluations of the FEM model. Alternatively, this extra computational cost can be reduced by using SQ. Only for the MC procedure it pays of\textcolor{black}{f} to use MOR, which is shown by comparing the times for the offline and online phase. The initial number of unknowns (8128) has been reduced to a basis of size~27. The difference in optimized volume for the various combinations of MOR and SQ is less then~$\unit{0.1}{\%}$. \textcolor{black}{The application of PSO and GA to the robust stochastic formulation results in smaller magnets, but has a computational cost that is more than 10 times higher than for Robust (i), even though the computations have been accelerated using parallelization. In contrast to SQP the PSO and GA algorithms do not make use of derivative information and thus evaluate much more machine models at every iteration step. Consequently, MOR is particularly beneficial in this case to speed up the computations, i.e., in the case of PSO the online costs are reduced by one order of magnitude, see Robust (iii) and (vii) in Table~\ref{tab:res_robopt}. All procedures terminated by reaching the desired accuracy.}

\textcolor{black}{To study the robustness of the optimized designs the failure rate has been determined. Around each optimum the same Monte Carlo sampling (with $N_\text{MC}=10000$) is performed with the same distribution as in \eqref{eq:distr}. For every sample the EMF is calculated and compared with $E_\mathrm{d}$. The failure rate is defined as $N_\text{fail}/N_{\text{MC}}$, where $N_\text{fail}$ is the number of machines with $E_0<E_\mathrm{d}$. Robust (i) gives} $\unit{4.36}{\%}$ of machines that do not fulfill $E_\mathrm{d}$, where for Linearized (i) model all machines fulfill the prerequisite. The linearized UQ setting differs only slightly from the full approach. \textcolor{black}{While yielding smaller magnets, the failure rates for the optima found by Robust (iii) and Robust (iv) are \textcolor{red}{marginally} smaller than for Robust (i) and Robust (ii).} One has to note that when the solution is linear in $\pp$, e.g. due to the linearization like in \eqref{eq:gen_opt_uq_rob_lin}, then $2\times2\times2$ collocation points are sufficient and no further approximation is introduced by the stochastic collocation. For the nonlinear solution such a coarse collocation grid also corresponds to a linearized model but it may differ from the one obtained by the Taylor expansion in~\eqref{eq:tay_exp}.

\begin{table*}[!t]
\processtable{Numerical results obtained for the nominal optimization procedures. The number in the column `Form' refers to the formulation in Section~\ref{sec:opt_prob}.\label{tab:res_nomopt}}
{\begin{tabular*}{\textwidth}{@{\extracolsep{\fill}}ccccccccccc@{}}\toprule
& ID& MOR& Method & UQ& Form &$\mathbf{P}$(mm)& Size (mm$^2$)& $E_0$ (V)& Failure rate (\%) & Time (s)\\
\midrule
Initial&- & - & -& - & -  &(19, 7, 7)& 133&  30.37& 50.17&-\\
\midrule
\multirow{5}{*}{Nominal}& (i) &- & SQP & -& \eqref{eq:opt} &(21.07, 2.98, 6.61)& 62.80& 30.37& 51.14& 3.3\\
&  (ii) &- & SQP& SQ& \eqref{eq:gen_opt_uq}& (21.07, 2.98, 6.61)& 62.80& 30.37& 50.19& 224\\
& (iii) &- &  SQP& MC& \eqref{eq:gen_opt_uq}&(21.07, 2.98, 6.61)& 62.80& 30.37& 50.19& 10774\\ 
& (iv)&- & \textcolor{black}{PSO}& -  & \eqref{eq:opt}&\textcolor{black}{(21.06, 2.98, 6.60)}& \textcolor{black}{62.80}& \textcolor{black}{30.37}& \textcolor{black}{51.34}& \textcolor{black}{818}\\
& (v)&- &\textcolor{black}{GA} & -  & \eqref{eq:opt}& \textcolor{black}{(21.44, 2.96, 7.16)}& \textcolor{black}{63.39}& \textcolor{black}{30.37}& \textcolor{black}{51.33}& \textcolor{black}{2385}\\
\botrule
\end{tabular*}}{}
\end{table*}
\begin{table*}[!t]
\processtable{Numerical results obtained for the robust optimization procedures with $\dbound=\unit{0.2}{mm}$. The column `MOR' \textcolor{red}{lists} the number of basis functions used, `Form' refers to the formulation in Section~\ref{sec:opt_prob} and `Time' distinguishes between online and offline costs if MOR was used. \label{tab:res_robopt}}
{\begin{tabular*}{\textwidth}{@{\extracolsep{\fill}}ccccccccccc@{}}\toprule
& ID& MOR& Method & UQ& Form & $\mathbf{P}$(mm)& Size (mm$^2$)& $E_0$ (V)& Failure rate (\%) & Time (s)\\
\midrule
Initial& -&  -& -&-  &- &(19, 7, 7)& 133& 30.37& 50.17&-\\
\midrule
\multirow{3}{*}{Linearized}& (i) &-& SQP  & -& \eqref{eq:opt_robust_1norm}& (20.88, 3.73, 6.82)& 77.86& 31.09& 0&  6.7\\
& (ii) &-& SQP& - & \eqref{eq:opt_robust4}& (20.95, 3.44, 6.78)&  72.01&  30.82& 4.28&  15\\
& (iii) & - & SQP& SQ& \eqref{eq:gen_opt_uq_rob_lin}& (20.97, 3.44, 6.82)&  72.10&  30.82& 4.03&  585\\
& (iv) & 27& SQP& SQ& \eqref{eq:gen_opt_uq_rob_lin}& (20.97, 3.44, 6.82)&  72.10&  30.82& 4.03&  234+301\\
\midrule
\multirow{8}{*}{Robust}& (i)&-& SQP & SQ &\eqref{eq:gen_opt_uq_rob}& (20.86, 3.53, 6.78)& 73.66& 30.82& 4.36&  239\\
&(ii) & -& SQP & MC & \eqref{eq:gen_opt_uq_rob}& (20.86, 3.53, 6.78)& 73.71& 30.81& 4.28& 14400\\
& (iii)&-&  \textcolor{black}{PSO} & SQ &\eqref{eq:gen_opt_uq_rob}& \textcolor{black}{(21.16, 3.22, 6.68)}& \textcolor{black}{68.04}& \textcolor{black}{30.86}& \textcolor{black}{3.91}&\textcolor{black}{2670}\\
& (iv)&-&  \textcolor{black}{GA} & SQ& \eqref{eq:gen_opt_uq_rob}& \textcolor{black}{(21.37, 3.23, 7.06)}& \textcolor{black}{69.02}& \textcolor{black}{30.85}&\textcolor{black}{4.02}& \textcolor{black}{3660}\\
& (v)& 27& SQP & SQ &\eqref{eq:gen_opt_uq_rob}& (20.86, 3.53, 6.78)& 73.66& 30.82& 4.44& 234+39\\
&(vi)& 27& SQP& MC &\eqref{eq:gen_opt_uq_rob}& (20.86, 3.53, 6.78)& 73.71& 30.82& 4.04& 234+1550\\
& (vii)& 27& \textcolor{black}{PSO} & SQ &\eqref{eq:gen_opt_uq_rob}& \textcolor{black}{(20.99, 3.27, 6.48)}& \textcolor{black}{68.52}& 30.84&\textcolor{black}{4.21}& \textcolor{black}{234+265}\\
& (viii) & 27&  \textcolor{black}{GA} & SQ&\eqref{eq:gen_opt_uq_rob}& \textcolor{black}{(21.60, 3.23, 7.39)}& \textcolor{black}{69.86}& 30.85&\textcolor{black}{4.17}& \textcolor{black}{234+1700}\\

\botrule
\end{tabular*}}{}
\end{table*}

\begin{figure}[!t]
\centering{\includegraphics[width=0.45\textwidth]{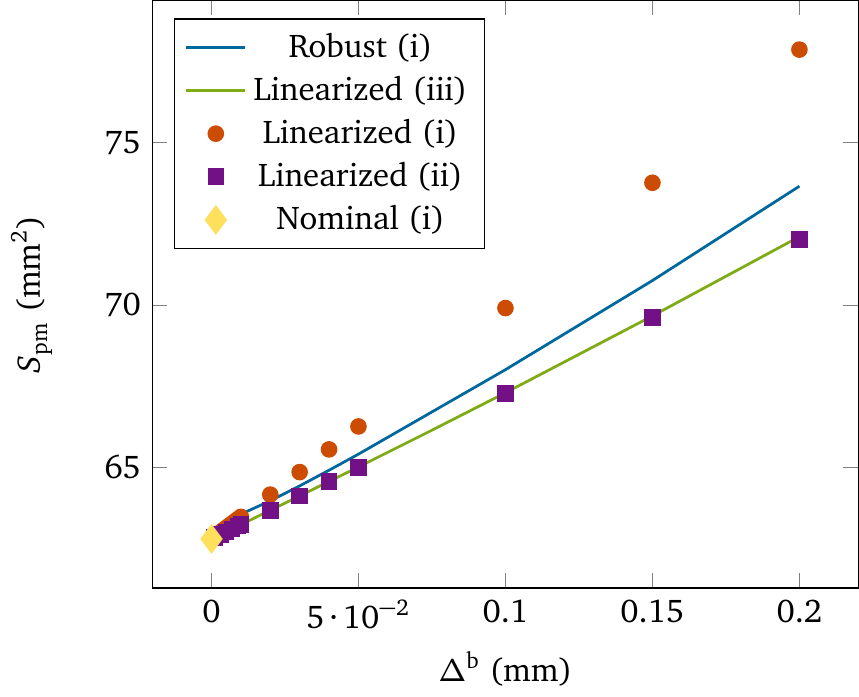}}
\caption{Optimization results for different values of $\dbound$. The yellow daimond is the size of the PMs obtained by using deterministic nominal optimization. The red circles and the purple squares are the results using deterministic robust optimization considering the 1-norm and the 2-norm respectively. The green line depicts the UQ optimization with the linearization. The blue line is the robust UQ optimization. \label{fig:delta_study1}}
\end{figure}


A visualization of the equivalence between the different \textcolor{black}{deterministic} approaches \textcolor{black}{for the PMSM example} is shown in Fig.~\ref{fig:delta_study1}. For this analysis the maximal deviation $\dbound$ is decreased to zero. This leads, as expected, to the optimized magnet size of Nominal~(i). The equivalence between Linearized~(ii) and Linearized~(iii)~can also be observed numerically, which \textcolor{black}{supports} Theorem 1. The results using robust optimization in the UQ setting using SQ and deterministic setting do differ. This is caused by the fact that Linearized~(i) is a more pessimistic scenario since it aims at mitigating the \textcolor{red}{worst-case}. By incorporating the second moment, one relies on more stochastic information during optimization which eventually translates into more optimistic results, because rare events are disregarded. This is also underlined by comparing the failure rates. 

\section{Conclusion \textcolor{red}{\& Outlook}}
The equivalence between robust worst-case optimization and variance-based optimization has been \textcolor{black}{mathematically derived} under the assumption that a linearization is applied and the norms are chosen adequately. Both approaches have been \textcolor{black}{numerically} compared using a simple benchmark problem, i.e., the reduction of the size of the permanent magnets in a permanent magnet synchronous machine while maintaining the electromotive force. It is found that robust optimization in the stochastic formulation gives less pessimistic results, since the \textcolor{red}{worst-case} might be unlikely to happen. However, the computational time is significantly increased whereas the implementation effort reduced since no (further) derivatives are needed. \textcolor{red}{The implementation of an affine decomposition facilitates the calculation of the derivatives and thus an efficient gradient based optimization procedure was obtained. The use of model order reduction has been shown to be beneficial when a lot of finite element evaluations are needed, which is the case for Monte Carlo sampling and when using the particle swarm optimization algorithm.}

\textcolor{red}{The robustness of the found optima was tested by determining the failure rates. It was found the worst-case optimization procedure ensured that all machines fulfill the prescribed quantity of interest. The other robust optimization approaches resulted in failrates around $\unit{4}{\%}$, whereas the nominal optimization formulations gave failrates around $\unit{50}{\%}$.}

\textcolor{red}{Even though the affine decomposition showed to be efficient, its application is not general e.g.~it cannot deal with rotation. Future work could focus on the use of design elements or iso-geometric analysis to construct efficient procedures to conduct shape optimization for electric machines with available gradients. Also a full 3D framework for the optimization procedure could be developed.}

\section*{Acknowledgment}
This work is supported by the German BMBF in the context of the SIMUROM project (grant no. 05M2013) and the PASIROM project (grant no. 05M18RDA), by the 'Excellence Initiative' of the German Federal and State Governments and the Graduate School of CE at TU Darmstadt.

\end{document}